\documentstyle[12pt]{article}

\def\qmod#1#2{{\hbox{}^{\displaystyle{#1}}}\!\big/\!\hbox{}_{
\displaystyle{#2}}}

\newfam\msbfam
\font\twelmsb=msbm10 at 12pt
\font\tenmsb=msbm10
\font\sevenmsb=msbm10 at 7pt
\font\fivemsb=msbm10 at 5pt

\textfont\msbfam=\twelmsb
\scriptfont\msbfam=\sevenmsb
\scriptscriptfont\msbfam=\fivemsb

\def\Bbb{\fam\msbfam\tenmsb}

\def\A{{\Bbb A}}

\def\C{{\Bbb C}}

\def\H{{\Bbb H}}

\def\L{{\Bbb L}}

\def\N{{\Bbb N}}

\def\P{{\Bbb P}}

\def\R{{\Bbb R}}
\def\Z{{\Bbb Z}}

\def\union{\mathop{\bigcup}}
\def\qed {\hfill\vrule height6pt width6pt depth0pt \bigskip}

\def\map{\longrightarrow}
\def\textmap#1{\mathop{\vbox{\ialign{
                                ##\crcr
    ${\scriptstyle\hfil\;\;#1\;\;\hfil}$\crcr
    \noalign{\kern-1pt\nointerlineskip}
    \rightarrowfill\crcr}}\;}}

\def\textlmap#1{\mathop{\vbox{\ialign{
                                ##\crcr
    ${\scriptstyle\hfil\;\;#1\;\;\hfil}$\crcr
    \noalign{\kern-1pt\nointerlineskip}
    \leftarrowfill\crcr}}\;}}

\newfam\meuffam
\font\twelmeuf=eufm10 at 12pt

\font\sevenmeuf=eufm7
\font\fivemeuf=eufm5

\textfont\meuffam=\twelmeuf
\scriptfont\meuffam=\sevenmeuf
\scriptscriptfont\meuffam=\fivemeuf

\def\germ{\fam\meuffam\twelmeuf}

\def\cg{{\germ c}}

\def\picture#1by#2(#3){
\vbox to #2 {
  \hrule width #1 height 0pt depth 0pt \vfill \special{picture #3}}
}

\def\scaledpicture#1by#2(#3scaled#4){{
\dimen0=#1  \dimen1=#2
\divide\dimen0 by 1000 \multiply\dimen0 by #4
\divide\dimen1 by 1000 \multiply\dimen1 by #4
\picture \dimen0 by \dimen1 (#3 scaled #4)}}
\def\dfigure#1by#2(#3scaled#4offset#5:#6)
  {\medskip
   \vglue 2mm minus 2mm
   $$
     \hbox{
       \hglue#5
       {\scaledpicture #1 by #2 (#3 scaled #4)}
     }
   $$
   \par\goodbreak
   \vglue 2mm minus 2mm
   \medskip}

\begin{document}
\def\Pr{{\rm Pr}}
\def\tr{{\rm Tr}}
\def\ad{{\rm ad}}
\def\End{{\rm End}}
\def\Pic{{\rm Pic}}
\def\NS{{\rm NS}}
\def\deg{{\rm deg}}
\def\Hom{{\rm Hom}}
\def\Aut{{\rm Aut}}
\def\Herm{{\rm Herm}}
\def\Vol{{\rm Vol}}
\def\pf{{\bf Proof: }}
\def\id{{\rm id}}
\def\im{{\rm im}}
\def\rk{{\rm rk}}
\def\coker{{\rm coker}}
\def\Sp{{\rm Sp}}
\def\Spin{{\rm Spin}}
\def\h{{\bf H}}
\def\dv{\bar\partial}
\def\ha{\frac{1}{2}}
\def\dva{\bar\partial_A}
\def\da{\partial_A}
\def\p{\partial\bar\partial}
\def\pa{\partial_A\bar\partial_A}
\def\Dr{{\raisebox{.16ex}{$\not$}}\hskip -0.35mm{D}}
\def\oo{{\scriptstyle{\cal O}}}
\def\ooo{{\scriptscriptstyle{\cal O}}}
\newtheorem{sz}{Satz}
\newtheorem{szfr}{Satzfr}
\newtheorem{thry}[sz]{Theorem}
\newtheorem{thfr}[szfr]{Th\'eor\`eme}
\newtheorem{pr}[sz]{Proposition}
\newtheorem{re}[sz]{Remark}
\newtheorem{co}[sz]{Corollary}
\newtheorem{cofr}[szfr]{Corollaire}
\newtheorem{dt}[sz]{Definition}
\newtheorem{lm}[sz]{Lemma}

\title{Seiberg-Witten invariants for manifolds with $b_+=1$, and the universal
wall crossing formula}

\author{Christian Okonek\thanks{Partially supported by: AGE-Algebraic
Geometry in Europe,
contract No ERBCHRXCT940557 (BBW 93.0187), and by  SNF, nr. 21-36111.92}
 \and     Andrei
Teleman$^*$\\ \\ Mathematisches Institut  Universit\"at
Z\"urich\\ Winterthurerstrasse 190, CH-8057 Z\"urich \\  \\ }
\date{  }
\maketitle
%

%

\section{Introduction}

The purpose of this paper  is to give a systematic description of the
Seiberg-Witten
invariants, which were introduced in [W], for manifolds with $b_+=1$. In
this situation,
compared to the general case $b_+>1$, several new features arise.

The Seiberg-Witten invariants for manifolds with $b_+>1$ are (non -
homogeneous) forms
$SW_{X,\oo}(\cg)\in \Lambda^* H^1(X,\Z)$, associated with an orientation
parameter
$\oo$ and a class of a $Spin^c(4)$-structures $\cg$ on $X$.

The invariants for manifolds with $b_+=1$ depend on a chamber structure;
they are
associated with  data $(\oo_1,{\bf H}_0,\cg)$, where $(\oo_1,{\bf H}_0)$ are
orientation parameters and $\cg$ is again the class of a
$Spin^c(4)$-structure on $X$. In
this case, the invariants are functions $SW_{X,(\oo_1,{\bf
H}_0)}(\cg):\{\pm\}\map
\Lambda^* H^1(X,\Z)$.

One of the main results of this paper is the proof of a universal wall
crossing formula.
This formula, which generalizes previous results of  [W], [KM]  and [LL]
%
%
describes the difference
$SW_{X,(\oo_1,{\bf H}_0)}(\cg)(+)-SW_{X,(\oo_1,{\bf H}_0)}(\cg)(-)$ as an
abelian
$Spin^c(4)$-form. More precisely, on elements $\lambda\in
\Lambda^r\left(\qmod{H_1(X,\Z)}{\rm Tors}\right)$ with $0\leq r \leq\min(
b_1,w_c)$,
we have:
$$\left[SW_{X,(\oo_1,{\bf H}_0)}(\cg)(+)-SW_{X,(\oo_1,{\bf
H}_0)}(\cg)(-)\right](\lambda)=\langle
\lambda\wedge\exp(-u_c),l_{\ooo_1}\rangle\ ,
$$
where $u_c\in \Lambda^2\left(\qmod{H_1(X,\Z)}{\rm Tors}\right)$ is given by
$u_c(a\wedge b)=\frac{1}{2}\langle a\cup b\cup c,[X]\rangle$, $a,b\in
H^1(X,\Z)$,
and $l_{\ooo_1}\in \Lambda^{b_1} H^1(X,\Z)$ represents the orientation
$\oo_1$ of
$H^1(X,\R)$. Here $c$ is the Chern class of $\cg$ and
$w_c:=\frac{1}{4}(c^2-3\sigma(X)-2e(X))$ is the index of $\cg$.

This formula has some important consequences, e.g. it shows that  Seiberg-Witten
invariants of manifolds with positive scalar curvature  metrics are essentially
topological invariants. According to Witten's vanishing theorem [W], one has
$SW_{X,(\oo_1,{\bf H}_0)}(\cg)(\pm)=0$ for at least one element of
$\{\pm\}$, and the other value is determined  be the wall crossing formula.

In the final part of the paper we show how to calculate $SW_{X,(\oo_1,{\bf
H}_0)}(\cg)(\pm)$ for K\"ahlerian surfaces. The relevant Seiberg-Witten
moduli spaces
have in this case a purely complex analytic description as Douady spaces
of curves
representing given homology classes: this description is essentially the
Kobaya\-shi-Hitchin correspondence obtained in [OT1].

Witten has shown that non-trivial invariants of K\"ahlerian surfaces with
$b_+>1$
must necessarily have index 0. This is not the case for surfaces with
$b_+=1$. We show
that a K\"ahlerian surface with $b_+=1$ and $b_1=0$ has $SW_{X,(\oo_1,{\bf
H}_0)}(\cg)(\{\pm\})=\{0,1\}$ or $SW_{X,(\oo_1,{\bf
H}_0)}(\cg)(\{\pm\})=\{0,-1\}$ as soon as the index of $\cg$ is
non-negative. For these
surfaces the invariants are therefore completely determined by their
reductions modulo
2.

There exist examples of 4-manifolds with $b_+=1$ which possess  - for every
prescribed
non-negative index - infinitely many classes $\cg$ of
$Spin^c(4)$-structures with
$SW_{X,(\oo_1,{\bf H}_0)}(\cg)\not\equiv 0$.


\section{ The twisted  Seiberg-Witten equations}

 Let $X$ be a closed
connected oriented    4-manifold, and let $c\in H^2(X,\Z)$ be a class with
$c\equiv
w_2(X)$ (mod 2). A
\underbar{compatible} $Spin^c(4)$-bundle is a    $Spin^c(4)$-bundle  $\hat
P$ over
$X$ with
$c_1(\hat P\times_{\det}\C)=c$ such that its $GL_+(4,\R)$-extension $\hat
P\times_{\tilde\pi}GL_+(4,\R)$ is isomorphic to the bundle of oriented frames in
$\Lambda^1_X$; here $\tilde\pi$ denotes the composition of the canonical
representation $\pi:Spin^c(4)\map SO(4)$ with the inclusion
$SO(4)\subset GL_+(4,\R)$. Let $\Sigma^{\pm}:=\hat P\times_{\sigma_{\pm}}\C^2$
be the associated spinor bundles with $\det\Sigma^{\pm}=\hat P\times_{\det}\C$
[OT1].
\begin{dt} A \underbar{Clifford}   \underbar{map}    of type $\hat P$ is
an orientation-preserving isomorphism  \hbox{ $\gamma:\Lambda^1_X\map
\hat P\times_{\pi}\R^4$}.
\end{dt}

The $SO(4)$-vector bundle $\hat P\times_\pi\R^4$ can be identified with the
bundle
$\R SU(\Sigma^+,\Sigma^-)$ of real multiples of $\C$-linear isometries  of
determinant 1  from $\Sigma^+$ to $\Sigma^-$.

 A  Clifford map
$\gamma$ defines a metric
$g_\gamma $ on $X$, a lift $\hat P\map P_{g_\gamma}$ of the associated
frame bundle,  and it induces isomorphisms
$\Gamma:\Lambda^2_{\pm}\map su(\Sigma^{\pm})$  [OT1].   We denote by ${\cal
C}={\cal C}(\hat P)$ the space of all Clifford maps of type $\hat P$. The
quotient
$\qmod{{\cal C}}{\im[\Aut(\hat P)\map \Aut(\hat P\times_\pi SO(4))]}$
para\-metrizes the set of all $Spin^c(4)$-structures of Chern class $c$,
whereas $\qmod{{\cal C}}{\Aut({\hat P\times_\pi SO(4))}}$ can be identified
with
the space ${\cal M}et_X$ of Riemannian metrics on $X$. In fact,
since ${\cal M}et_X$ is contractible, we have  a  natural  isomorphism
$\qmod{{\cal C}}{ \Aut( \hat P)}
\textmap{\simeq}{\cal M}et_X\times
\pi_0\left(\qmod{{\cal C}} { \Aut( \hat P)} \right)
$, where the second factor is a ${\rm Tors}_2H^2(X,\Z)$-torsor; it
para\-metrizes the set of equivalence classes of $Spin^c(4)$ -structures with
Chern class $c$ on $(X,g)$,  for an arbitrary  metric $g$. The latter
assertion follows
from the fact that the map $H^1(X, {\Z}_{2})\map
H^1(X,\underline{Spin^c}(4))$ is
trivial for any 4-manifold $X$.
We use the symbol $\cg$ to denote elements in $\pi_0\left(\qmod{{\cal C}}{ \Aut(
\hat P)}\right)$, and we denote by $\cg_\gamma$ the connected component defined
by
$[\gamma]\in\qmod{{\cal C}} { \Aut( \hat P)}$.

A fixed Clifford map $\gamma$ defines a bijection
between unitary connections in
$\hat P\times_{\det}\C$ and $Spin^c(4)$-connections in $\hat P$ which lift (via
$\gamma$) the Levi-Civita connection in
$P_{g_\gamma}$, and allows to associate  a Dirac operator $\Dr_A$ to  a
connection
$A\in{\cal A}(\hat P\times_{\det}\C)$.
\begin{dt} Let $\gamma$ be a Clifford map, and let $\beta\in Z^2_{\rm
DR}(X)$ be a
closed 2-form. The
$\beta$-twisted Seiberg-Witten equations are
$$\left\{\begin{array}{lll}
\Dr_{A}\Psi&=&0\\
\Gamma\left((F_A+{2\pi i}\beta)^+\right)&=&2(\Psi\bar\Psi)_0\ .
\end{array}\right. \eqno{(SW^{\gamma}_\beta)}
$$
\end{dt}
These twisted Seiberg-Witten equations arise
 naturally in
connection with certain non-abelian monopoles [OT2], [T].  They should
\underbar{not} be regarded as perturbation of $(SW^\gamma_0)$, since later the
cohomology class of $\beta$ will be fixed.

Let ${\cal W}_{X,\beta}^{\gamma} $  be the moduli space  of
solutions $(A,\Psi)\in{\cal A}(\det \Sigma^+)\times A^0(\Sigma^+)$ of
$(SW^{\gamma}_\beta)$ modulo  the natural action
$((A,\Psi),f)\longmapsto
(A^{f^2},f^{-1}\Psi)$ of the gauge group
${\cal G}={\cal C}^{\infty}(X,S^1)$.

Since two Clifford maps lifting the same pair $(g,\cg)$ are equivalent modulo
$\Aut(\hat P)$, the moduli space  ${\cal W}_{X,\beta}^{\gamma}$ depends up to
\underbar{canonical} isomorphism only on  $(g_\gamma,\cg_\gamma)$ and $\beta$.

Now fix a class $b\in H^2_{\rm DR}(X)$, consider
$(SW^{\gamma}_\beta)$ as equation for a triple $(A,\Psi,\beta)\in {\cal
A}(\det \Sigma^+)\times A^0(\Sigma^+)\times b$, and let
${\cal W}_{X,b}^{\gamma}\subset\qmod{ {\cal A}(\det \Sigma^+)\times
A^0(\Sigma^+) \times b }{{\cal G}}$
be the (infinite dimensional) moduli space of solutions.  Finally we need the
universal moduli space
${\cal W}_X \subset\qmod{ {\cal A}(\det \Sigma^+)\times
A^0\ (\Sigma^+)\times Z^2_{DR}(X)\times {\cal C} }
{{\cal G}}$
of solutions of $(SW^{\gamma}_\beta)$ regarded as equations for tuples
$(A,\Psi,\beta,\gamma)\in{\cal A}(\det \Sigma^+)\times
A^0(\Sigma^+)\times Z^2_{DR}(X)\times {\cal C} $.

We complete the spaces  ${\cal A}(\det\Sigma^+)$, $A^0(\Sigma^{\pm})$ and $A^2$
with respect to the Sobolev norms $L^2_q$, $L^2_q$ and $L^2_{q-1}$, and
the gauge group ${\cal G}$ with respect to $L^2_{q+1}$,
but we suppress the Sobolev subscripts in our notations. As usual we denote
by the
superscript ${\ }^*$ the open subspace of a moduli space where the spinor
component is non-zero.
\begin{dt} Let $c\in H^2(X,\Z)$ be   characteristic. A pair
$(g,b)\in{\cal M}et_X\times H^2_{\rm DR}(X)$
 is  $c$-good if  the g-harmonic
representant  of $(c-b)$ is not $g$-anti-selfdual.
\end{dt}
A pair  $(g,b)$ is $c$-good for every metric $g$ if
$(c-b)\ne\ 0$ and $(c-b)^2\geq 0$.

\begin{pr} Let $X$ be a closed  oriented 4-manifold, and let $c\in H^2(X,\Z)$ be
characteristic. Choose a compatible $Spin^c(4)$-bundle $\hat P$  and  an
element $\cg\in \pi_0\left(\qmod{ {\cal C}}{\Aut(
\hat P)} \right)$. \\
i)  The projections
$p:{\cal W}_X \map Z^2_{DR}(X)\times{\cal C}$ and   $p_{\gamma,b}:{\cal
W}_{X,b}^{\gamma } \map b$  are proper for all choices of $\gamma$ and $b$.\\
ii)     ${\cal W}_X^* $  and ${{\cal W}_{X,b}^{\gamma}}^* $  are  smooth
manifolds  for all $\gamma$ and $b$ .\\
iii)    ${{\cal W}_{X,b}^{\gamma}}^* ={\cal W}_{X,b}^{\gamma} $ if
$(g_\gamma,b)$ is
$c$-good.\\
iv)  If $(g_\gamma,b)$ is $c$-good, then every pair $(\beta_0,\beta_1)$ of
regular
values of $p_{\gamma,b}$ can be joined by a smooth path
$\beta:[0,1]\map b$   such that the fiber product
$[0,1]\times_{(\beta, p_{\gamma,b})}{\cal W}_{X,b}^{\gamma} $ defines a
smooth  bordism between   ${\cal W}_{X,\beta_0}^{\gamma} $ and ${\cal
W}_{X,\beta_1}^{\gamma} $. \\
v) If   $(g_0,b_0)$,  $(g_1,b_1)$ are $c$-good pairs which can be joined by a
smooth path  of $c$-good  pairs,   then there is a
smooth path $(\beta,\gamma):[0,1]\map   Z^2_{\rm DR}(X)\times{\cal
C}$ with the following properties:\\
\hspace*{0.8cm}1. $[\beta_i]=b_i$ and $g_{\gamma_i}=g_i$ for $i=0,  1$.  \\
\hspace*{0.8cm}2. $\gamma_t$ lifts $(g_{\gamma_t},\cg)$ and
$(g_{\gamma_t},[\beta_t])$ is $c$-good for every $t\in[0,1]$.\\
\hspace*{0.8cm}3.   $[0,1]\times_{((\beta,\gamma),p)}{\cal
W}_X ^*$ is a smooth bordism between
${\cal W}_{X,\beta_0}^{\gamma_0} $ and
${\cal W}_{X,\beta_1}^{\gamma_1} $.\\
vi) If $b_+>1$, then any two $c$-good pairs $(g_0,b_0)$,
$(g_1,b_1)$  can be joined by a smooth path of $c$-good pairs .
\end{pr}
\pf \\
{\sl i)} See  [KM], Corollary 3.\\
{\sl ii)} It suffices to show that, for a fixed class $b\in H^2_{\rm
DR}(X)$, the map
$$F:{\cal A}(\det\Sigma^+)\times A^0(\Sigma^+)\times b\map
A^0(\Sigma^-)\times  iA^0(su(\Sigma^+))
$$
given by
$F(A,\Psi,\beta)=\left(\Dr_A\Psi,\Gamma\left((F_A+2\pi
i\beta)^+\right)-2(\Psi\bar\Psi)_0\right)
$
is a submersion in every point $\tau=(A,\Psi,\beta)$ with $\Psi\ne 0$ and
$\Dr_{A}\Psi=0$. To see this write $F=(F^1,F^2)$ for the components of
$F$, and consider a pair
$(\Psi^-,S)\in A^0(\Sigma^-)\times iA^0(su(\Sigma^+)$ which is
$L^2$-orthogonal to the image of $d_\tau(F)$. Using variations of $\beta$ by
exact forms, we see that $\Gamma^{-1}(S)\in i A^2_+$ is orthogonal to
$id^+(A^1)$, hence  must be a  harmonic selfdual form. This implies
$\left\langle d_\tau(F^2)(a,0,0),S\right\rangle=0$ for any variation
$a\in iA^1$, and therefore
$$Re\langle \gamma(a)(\Psi),\Psi^-\rangle=  Re\left\langle
d_\tau(F^1)(a,0),(\Psi^-)\right\rangle=\left\langle
d_\tau(F)(a,0,0),(\Psi^-,S)\right\rangle=0
$$
for all $a\in iA^1$.   Since $\Psi$ is
Dirac-harmonic and non-trivial, it cannot vanish on   non-empty open sets.
Therefore,  the multiplication map $\gamma(\cdot)\Psi:iA^1\map
A^0(\Sigma^-)$ has
$L^2$-dense image, so that we must have $\Psi^-=0$. Using now the vanishing of
$\Psi^-=0$, we get
$$\langle d_\tau F^2(0,\psi),S\rangle=\langle
d_\tau F(0,\psi,0),(\Psi^-,S)\rangle=0
$$
for all variations $\psi\in A^0(\Sigma^+)$. Since $d_\tau
F^2(0,\cdot):A^0(\Sigma^+)\map  i A^0(su(\Sigma^+))$ has $L^2$-dense image,
we must also  have $S=0$.\\
%
{\sl iii)} This is clear since  the  form $F_A+2\pi i\beta$ represents the
cohomology class
$-2\pi i(c-b)$: If $(A,0)$ was a solution of $(SW^{\gamma}_\beta)$, then
the $g_\gamma$-anti-selfdual form $\frac{i}{2\pi}F_A-\beta$ would be the
$g_\gamma$-harmonic representant of $c-b$. \\
{\sl iv)} This  follows from {\sl ii)}, {\sl iii)}  and Proposition
(4.3.10) of  [DK].\\
{\sl v)} This is a consequence of {\sl ii)}, {\sl iii)}, Proposition
(4.3.10) of  [DK], and
the  fact that the condition "$c$-good" is open for fixed $c$.\\
{\sl vi)}  For fixed  $c$,   the  closed subspace of $Z^2_{\rm
DR}(X)\times{\cal C}$
consisting of pairs  $(\beta, {\gamma})$ with $(g_\gamma,[\beta])$ not $c$-good
has codimension
$b_+$. Its complement is therefore  connected when $b_+\geq 2$. This fact
was  noticed by C. Taubes [Ta].
\qed
\\
{\bf Remark:} Suppose $(g,b)$ is a $c$-good pair. Then the bordism type of
a smooth
moduli space ${\cal W}_{X,\beta}^{\gamma}$ with $g_\gamma=g$ depends only
on the pair $(g,\cg)$ and the class $b$ of $\beta$. Furthermore, this
bordism type
does not change as long as one varies $(g,b)$ in a smooth 1-parameter family of
$c$-good pairs. Note that the statement makes sense since the set
$\pi_0(\qmod{{\cal C}}{\Aut(\hat P)})$ to which $\cg$ belongs was
defined independently of the metric.
\\

 \section{ Seiberg-Witten invariants for 4-manifolds with $b_+=1$}
 Let $X$
be a closed connected   oriented 4-manifold , $c$ a characteristic element,  and
$\hat P$ a compatible $Spin^c(4)$-bundle. We put
$${\cal B}(c)^*:= \qmod{{\cal A}(\det\Sigma^+)
\times(A^0(\Sigma^+)\setminus\{0\})}{{\cal G}}\ .$$
\begin{lm}
 ${\cal B}(c)^*$ has the weak homotopy type of $K(\Z,2)\times
K(H^1(X,\Z),1)$. There is a natural isomorphism
$$ \nu:\Z[u]\otimes \Lambda^*(\qmod{H_1(X,\Z)}{\rm Tors})\map  H^*({\cal
B}(c)^*,\Z)
.$$
\end{lm}
\pf The inclusion $S^1\subset{\cal G}$ defines a canonical exact sequence
$$1\map S^1\map {\cal G}\map\overline{{\cal G}}\map 1
$$
with $\overline{{\cal G}}:=\qmod{{\cal G}}{ S^1}$, and the exponential map
yields
a natural identification of $\overline{{\cal G}}$ with the product
$\qmod{{\cal C}^{\infty}(X,\R)}{\R}\times H^1(X,\Z)$.  The choice of a base
point
$x_0\in X$ induces a splitting $ev_{x_0}:{\cal G}\map  S^1$ of the exact
sequence,
and therefore a homotopy equivalence of classifying spaces $B{\cal G}\map
BS^1\times B\overline{{\cal G}}$; the homotopy class of this map is
independent of
$x_0$ when $X$ is connected. Since  ${\cal
A}(\det \Sigma^+)\times(A^0(\Sigma^+)\setminus\{0\})$ is weakly contractible,
${\cal B}(c)^*$ has the weak homotopy type of $B{\cal G}$. We fix weak homotopy
equivalences ${\cal B}(c)^*\simeq B{\cal G}$,
$BS^1\simeq K(\Z,2)$ and $B\overline{{\cal G}}\simeq K(H^1(X,\Z),1)$ in the
natural homotopy classes. Since the homotopy class  of the induced weak
homotopy equivalence ${\cal B}(c)^*\map K(\Z,2)\times K(H^1(X,\Z),1)$ is canonical,
we obtain a natural isomorphism
$$H^*({\cal
B}(c)^*,\Z)\simeq H^*( K(\Z,2))\otimes
H^*(K(H^1(X,\Z),1))\simeq \Z[u]\otimes\Lambda^*(\qmod{H_1(X,\Z)}{\rm Tors})\ .$$
 \qed

{\bf Remark:}  Let ${\cal G}_0$ be the kernel of the evaluation map
$ev_{x_0}:{\cal
G}\map  S^1$, and set ${\cal B}_0(c)^*:=\qmod{{\cal A}(\det\Sigma^+)
\times(A^0(\Sigma^+)\setminus\{0\})}{{\cal G}_0}$. The group $\qmod{{\cal
G}}{{\cal G}_0}\simeq S^1$ acts freely on  ${\cal B}_0(c)^*$ and defines a
principal
$S^1$-bundle over  ${\cal B}(c)^*$. The first Chern class of this bundle is
the class
$u$ defined above.\\

Suppose now that $(g,b)$ is a c-good pair, and  fix $\cg\in\pi_0(\qmod{{\cal
C}}{\Aut(\hat P)})$ . The moduli space
${\cal W}^{\gamma}_{X,\beta} $ is  a compact manifold of dimension
$w_c:=\frac{1}{4}(c^2-2e(X)-3\sigma(X))$  for every  lift $\gamma$ of $(g,\cg)$
and every regular value
$\beta$ of $p_{\gamma,b}$. It can be oriented by using
the canonical complex orientation of the line bundle
$\det{}_\R({\rm index} (\Dr)) $
over
${\cal B}(c)^*$ together with a chosen orientation $\oo$ of the line
$\det(H^1(X,\R))\otimes\det(\H^2_{+,g}(X)^{\vee})$. Let $[{\cal
W}^{\gamma}_{X,\beta}  ]_{{\raisebox{-.5ex}{$\oo$}}}\in H_{w_c}({\cal
B}(c)^*,\Z)$ be the fundamental class associated  with the choice of $\oo$.

\begin{dt} Let $X$ be a  closed connected oriented 4-manifold, $c\in
H^2(X,\Z)$ a
characteristic element, $(g,b)$ a $c$-good pair, $\cg\in\pi_0(\qmod{{\cal
C}}{\Aut(\hat P)})$, and $\oo$ an orientation of the line
$\det(H^1(X,\R))\otimes\det(\H^2_{+,g}(X)^{\vee})$. The corresponding
Seiberg-Witten form is the element
$$SW_{X,\oo}^{(g,b)}(\cg)\in \Lambda^* H^1(X,\Z)
$$
defined by
$$SW_{X,\oo}^{(g,b)}(\cg)(l_1\wedge\dots\wedge l_r)=\left\langle\nu(l_1)
\cup\dots\cup\nu(l_r)\cup u^{\frac{w_c-r}{2}}, [{\cal W}^{\gamma}_{X,\beta}
]_{{\raisebox{-.5ex}{$\oo$}}} \right\rangle\
$$
for decomposable elements
$l_1\wedge\dots\wedge l_r$ with  $r\equiv w_c$ (mod 2) .
Here $\gamma$ lifts the pair $(g,\cg)$ and $\beta\in b$ is a regular value of
$p_{\gamma,b}$.
\end{dt}
{\bf Remark:} The form $SW_{X,\oo}^{(g,b)}(\cg)$ is well-defined, since the
 cohomology classes $u$, $\nu(l_i)$, as well as the trivialization of the
orientation line bundle extend  to   the  quotient $\qmod{{\cal A}(\det
\Sigma^+)\times (A^0(\Sigma^+)\setminus\{0\})\times b}{{\cal G}}$,  and
since, by
Proposition 4 {\sl iv)}, the homology class defined by $[{\cal
W}^{\gamma}_{X,\beta} ]_{{\raisebox{-.5ex}{$\oo$}}}$ in this
quotient depends only on $(g_\gamma,\cg_\gamma)$ and $b$.\\

Now there are two cases:

If $b_+>1$, then,  by Proposition 4 {\sl v)},  {\sl vi)},  the
form $SW_{X,\oo}^{(g,b)}(\cg)$ is also independent  of $(g,b)$,  since the
cohomology classes and the trivialization of the orientation line bundle
extend to
$\Aut(\hat P)$-invariant objects on the universal quotient $\qmod{{\cal
A}(\det\Sigma^+)\times (A^0(\Sigma^+)\setminus\{0\})\times Z^2_{\rm
DR}(X)\times{\cal C}}{{\cal G}}$. Thus we may simply write
$SW_{X,\oo}(\cg)\in
\Lambda^* H^1(X,\Z)$. If $b_1=0$, then we obtain numbers which we denote by
$n_{\cg}^{\ooo}$; these numbers can be considered as refinements of the numbers
$n_c^{\ooo}$ which were defined in [W]. Indeed, $n_c^{\ooo}=\sum_{\cg}
n^{\ooo}_{\cg}$, the summation being over all $\cg\in\pi_0(\qmod{{\cal
C}}{\Aut(\hat P)})$.

Suppose now that $b_+=1$. There is a natural map ${\cal M}et_X\map
\P(H^2_{\rm DR}(X))$ which sends a metric $g$ to the line $\R[\omega_+]\subset
H^2_{\rm DR}(X)$, where $\omega_+$ is any non-trivial $g$-selfdual harmonic
form. Let ${\bf H}$ be the hyperbolic space
$${\bf H}:=\{\omega\in H^2_{\rm DR}(X)|\ \omega^2=1\}\ .
$$

${\bf H}$ has two connected components, and the choice of one of them
orients the
lines $\H^2_{+,g}(X)$ for all metrics $g$. Furthermore, having fixed a component
${\bf H}_0$ of ${\bf H}$, every metric defines a unique $g$-self-dual form
$\omega_g$ with $[\omega_g{]\in\bf H}_0$.
\begin{dt} Let $X$ be a manifold with $b_+=1$, and let $c\in H^2(X,\Z)$ be
characteristic. The \underbar{wall} associated with $c$ is the hypersurface
$c^{\bot}:=\{(\omega,b)\in{\bf H}\times H^2_{\rm DR}(X)|\
(c-b)\cdot\omega=0\}$.
The connected components of ${\bf H}\setminus c^{\bot}$ are called
\underbar{chambers} of type $c$.
\end{dt}

Notice that the walls are non-linear! Every characteristic element $c$ defines
precisely four chambers of type $c$, namely
$$C_{{\bf H}_0,\pm}:=\{(\omega,b)\in{\bf H}_0\times H^2_{\rm DR}(X)|\
\pm(c-b)\cdot\omega<0\}\ ,
$$
 where ${\bf H}_0$ is one of the components of ${\bf H}$. Each of   these four
chambers contains pairs of the form $([\omega_g],b)$. Let $\oo_1$ be an
orientation of $H^1(X,\R)$.
\begin{dt} The Seiberg-Witten invariant associated with the data $(\oo_1,{\bf
H}_0,\cg)$ is the function
$$
SW_{X,(\oo_1,{\bf H}_0)}(\cg):\{\pm\}\map  \Lambda^* H^1(X,\Z)
$$
given by  $SW_{X,(\oo_1,{\bf H}_0)}(\cg)(\pm):=SW^{(g,b)}_{X,\oo}(\cg)$, where
$\oo$ is the orientation defined by $(\oo_1,{\bf H}_0)$, and $(g,b)$ is a
pair  such
that $([\omega_g],b)$  belongs to the chamber $C_{{\bf H}_0,\pm}$.
\end{dt}
{\bf Remark: }  The intersection $c^{\bot}\cap {\bf H}\times\{0\}$ defines a
non-trivial wall in ${\bf H}\times\{0\}$ if and only if $c^2<0$. This means that
invariants of index $w_c<\frac{b_2-10}{4}+b_1$ could also be defined for the
chambers $C^{\pm}_{{\bf H}_0}:=\{\omega\in{\bf H}_0|\ \pm c\cdot \omega<0\}$.
However, when $c$ is rationally non-zero and $c^2\geq 0$, then ${\bf
H}_0\times\{0\}$ is entirely contained in one of the chambers $C_{{\bf
H}_0,\pm}$.\\

Note that, changing the orientation $\oo_1$  changes  the invariant  by a factor
$-1$, and that $SW_{X,(\oo_1,-{\bf H}_0)}(\cg)(\pm)=-SW_{X,(\oo_1, {\bf
H}_0)}(\cg)(\mp)$.\\ \\
{\bf Remark: } A different approach - adapting ideas from intersection theory to
construct "Seiberg-Witten multiplicities" - has been proposed by R. Brussee.\\

\section{The wall crossing formula}
In this section we prove a wall crossing formula  for the complete
Seiberg-Witten
invariant $SW_{X,(\oo_1,{\bf H}_0)}(\cg)$. This formula generalizes previous
results of [W],  [KM] and [LL].  Our proof is based on ideas similar to the
ones in  [LL],
but our method - using the real bow up
of the (singular)  locus of reducible points as in  [OT2] - has some advantages:
It allows us to construct explicitely a smooth bordism to which the cohomology
classes $u$, $\nu(l_i)$ extend in  a natural way, and it enables us a  to
give  a simple description of that part of  its boundary which lies  on the
wall.

Our main result can be  formulated by saying that the the difference
$SW_{X,(\oo_1,{\bf H}_0)}(\cg)(+)-SW_{X,(\oo_1, {\bf H}_0)}(\cg)(-)$ is an
abelian $Spin^c$-form. These abelian $Spin^c$-forms  are topological invariants
which will be introduced in the first subsection. Using a remark of [LL],
we give an
explicit formula for these invariants in the case of manifolds with $b_+=1$.
\subsection{Abelian $Spin^c(4)$-forms }

Let $X$ be an oriented 4-manifold, $c\in H^2(X,\Z)$ a characteristic
element, and
$\hat P$ a compatible
$Spin^c(4)$-bundle. Let $\Sigma^{\pm}$ be the associated spinor bundles, define
$L:=\det\Sigma^{\pm}$, and put:
$${\cal B}(L):=\qmod{{\cal A}(L)}{{\cal G}_0}\  ,\ \  {\cal
B}'(L):=\qmod{{\cal A}(L)}{{\cal G}_0^2}\ .$$
There is an obvious covering projection
$$s:{\cal B}'(L)\map {\cal B}(L)$$
with fiber \ $\qmod{{\cal G}_0}{{\cal G}_0^2}=\qmod{H^1(X,\Z)}{2H^1(X,\Z)}$.

The homotopy equivalence ${\cal G}_0\textmap{\simeq} H^1(X,\Z)$  induces
canonical
isomorphisms
$$\mu:\Lambda^*\left(\qmod{H_1(X,\Z)}{\rm
Tors})\right)\textmap{\simeq}  H^*({\cal B}(L),\Z)\ ,$$
$$\nu':\Lambda^*\left(\qmod{H_1(X,\Z)}{\rm
Tors})\right)\textmap{\simeq}  H^*({\cal B}'(L),\Z)\ ,$$
such that $s^*\circ\mu=2\nu'$.

Now fix a metric $g$ and let $h_c$ be the harmonic representant of the de Rham
 class $c_{\rm DR}\in H^2_{\rm DR}(X)$. The spaces
${\cal B}(L)$, ${\cal B}'(L)$ are  trivial  fibre bundles over the affine
subspace
$h_c^+ +d^+(A^1)=h_c^++(\H^2_{g,+}(X))^{\bot}\subset A^2_+$ via the maps induced by
$A \longmapsto \frac{i}{2\pi}F_A^+$.

For a given 2-form $\beta\in h_c^+ +d^+(A^1)$  let
$${\cal T}_\beta(L)\subset  {\cal B}(L)\ ,\ \ {\cal
T}_{\beta}'(L)\subset{\cal B}'(L)$$
 be the fibers of these maps over $\beta$.  These fibers are tori, consisting of
equivalence classes of solutions of   of the equation
$$F_A^++2\pi i\beta=0
$$
modulo ${\cal G}_0$ respectively ${\cal G}_0^2$.   Indeed, the choice of a
solution $A_0\in{\cal A}(L)$   yields identifications
$${\cal T}_\beta(L)=\qmod{H^1(X,\R)}{H^1(X,\Z)}\ ,\ \ {\cal
T}_{\beta}'(L)=\qmod{H^1(X,\R)}{2H^1(X,\Z)}\ .$$

Now fix a Clifford map $\gamma\in{\cal C}$ of type $\hat P$ with $g_\gamma=g$,
and a 2-form $\beta\in h_c^+ +d^+(A^1)$. Let $S(A^0(\Sigma^+))$ be the
unit  sphere in $A^0(\Sigma^+)$ with respect to the $L^2$-norm.
\begin{dt} The $\beta$-twisted $Spin^c(4)$-equations for a pair
$(A,\Phi)\in {\cal
A}(L)\times S(A^0(\Sigma^+))$ are the equations
$$\left\{\begin{array}{ccl}\Dr_A \Phi&=&0\\
 F_A^++2\pi i\beta&=&0 \ .
\end{array}\right.  \eqno{(S^\gamma_\beta)}
$$
\end{dt}
 The gauge group ${\cal G}$ acts on ${\cal A}(L)\times S(A^0(\Sigma^+))$  by
$(A,\Phi)\cdot f= (A^{f^2},f^{-1}\Phi)$, letting invariant the space of
solutions of $(S^\gamma_\beta)$. We denote by $\hat {\cal T}_{\beta}'^\gamma(L)$
the moduli space of solutions; it is    a
projective fiber space over the "Brill-Noether locus" in ${\cal
T}_{\beta}'(L)$.
Let
$$q:\qmod{{\cal A}(L)\times S(A^0(\Sigma^+))}{{\cal G}}\map  {\cal B}'(L)
$$
be the natural projection.
\begin{lm} There is a natural isomorphism
$$\Z[u]\otimes
\Lambda^* \left(\qmod{H_1(X,\Z)}{\rm Tors}\right)\map H^*\left(\qmod{{\cal
A}(L)\times S(A^0(\Sigma^+))}{{\cal G}},\Z\right)
$$
whose restriction to $\Lambda^* \left(\qmod{H_1(X,\Z)}{\rm Tors}\right)$ factors
as $q^*\circ\nu'$.   The class $u$ restricts to the positive generator of the
 second cohomology group $\P(A^0(\Sigma^+))$ of the fibers of $q$.
\end{lm}
This lemma, as well as the following proposition can be proved using the same
methods as in the proofs of Lamma 5 and Proposition 4.

Set $\delta_c:=\frac{1}{8}(c^2-\sigma(X))$.
\begin{pr}  For generic elements $\beta\in h_c^+ + d^+(A^1)$, the moduli
space $\hat {\cal
T}_{\beta}'^\gamma(L)$ is a  closed smooth manifold of dimension
$b_1+2\delta_c-2$. It can be oriented by choosing an orientation
$\oo_1$ of $H^1(X,\R)$. The fundamental class
$$[\hat {\cal
T}_{\beta}'^\gamma(L)]_{\ooo_1}\in
H_{b_1+2\delta_c-2}\left( \qmod{{\cal A}(L)\times S(A^0(\Sigma^+))}{{\cal
G}} ,\Z\right)$$
 depends only on the component $\cg_\gamma\in\pi_0\left(\qmod{{\cal
C}}{\Aut(\hat P)}\right)$.
\end{pr}
We can now define the abelian $Spin^c(4)$-forms.
\begin{dt}  Let $X$ be a closed connected oriented 4-manifold, and $\oo_1$ an
orientation of $H^1(X,\R)$. Let $c\in H^2(X,\Z)$ be a characteristic element and
$\cg\in\pi_0\left(\qmod{{\cal C}}{\Aut(\hat P)}\right)$.  The corresponding
$Spin^c(4)$-form is the element $\sigma_{X,{\ooo_1}}(\cg)\in
\Lambda^*(H^1(X,\Z))$ defined  by the formula
$$\sigma_{X,\ooo_1}(\cg)(l_1\wedge\dots\wedge
l_r):=\langle\nu'(l_1)\cup\dots\cup\nu'(l_r)\cup
u^{\frac{b_1+2\delta_c-2-r}{2}},[\hat {\cal
T}_{\beta}'^\gamma(L)]_{\ooo_1}\rangle
$$
for   decomposable elements $l_1\wedge\dots\wedge l_r$  with $r\equiv b_1$
(mod 2).
Here $\gamma$ induces $\cg=\cg_\gamma$ and $\beta$ is generic.
\end{dt}

Note that the expected dimension $b_1+2\delta_c-2$ of the moduli space
  $\hat {\cal T}_{\beta}'^\gamma(L)$ coincides with $w_c$ iff
$b_+=1$. The $Spin^c$-forms $\sigma_{X,\ooo_1}(\cg)$ are topological invariants,
which can be explicitely computed.  This is our next aim.\\

Let $pr_2:{\cal A}(L)\times X\map  X$ be the projection onto the second
factor, and put
$$\P:=\qmod{pr_2^*(\hat P)}{{\cal G}_0}\ ,
$$
where ${\cal G}_0$ acts on $pr_2^*(\hat P)={\cal A}(L)\times\hat P$ by
  $(A,\hat p)\cdot f =(A^{f^2}, f(\hat p))$. The universal
$Spin^c(4)$-bundle  $\P$  over ${\cal B'}(L)\times X$  comes with a tautological
connection $\A$ in the
$X$-direction.     Let $\L:=\qmod{pr_2^*(L)}{{\cal G}_0}$ be the universal line
bundle over ${\cal B}(L)\times X$; its pull back  $(s\times\id)^*(\L)$ to ${\cal
B}'(L)\times X$ is the universal determinant bundle $\P\times_{\det}\C$.  The
Chern class  $c_1(\L)$ has a K\"unneth decomposition $c_1(\L)=1\otimes c+
c_1(\L)^{1,1}$ whose $(1,1)$-component $c_1(\L)^{1,1}\in H^1({\cal
B}(L),\Z)\otimes H^1(X,\Z)$, considered as homomorphism
$\mu_1:\qmod{H^1(X,\Z)}{\rm Tors}\map  H^1({\cal B}(L),\Z)$, is given by the
restriction of the isomorphism $\mu$.

We fix a basis $(l_i)_{1\leq i\leq b_1}$ of
$\qmod{H_1(X,\Z)}{\rm Tors}$  and let $l^i$ be the elements of the dual
basis. Then
$$c_1(\L)=1\otimes c+\sum\limits_{i=1}^{b_1} \mu(l_i)\otimes l^i\ ,$$
$$c_1(\P\times_{\det}\C)=1\otimes c+2\sum\limits_{i=1}^{b_1}\nu'(l_i)\otimes
l^i  \ .$$

Let $\gamma:\Lambda^1\map  \hat P\times_{\pi}\R^4$ be a Clifford map  and
 $\raisebox{-.18ex}{$\gamma$}\hskip-5.5pt{\gamma}:pr_2^*(\Lambda^1)\map
\P\times_{\pi}\R^4$ the induced isomorphism.
We  restrict the universal objects
$(\P,\raisebox{-.18ex}{$\gamma$}\hskip-5.5pt{\gamma},\A)$ to the subspace
${\cal T}'_{\beta}(L)\times X$ and   denote the restrictions by the same
symbols.   The Chern character of the  virtual bundle
$index(\Dr_\A)$ over the torus ${\cal T}'_{\beta}(L)$  is
$$ch(index(\Dr_\A))=\left[\left(e^{\frac{c}{2}+\sum \nu'(l_i)\otimes l^i
}\right)
\cup \left(1-\frac{1}{24}p_1(X)\right)\right]/[X]=
$$
$$=\left(e^{\frac{c}{2}+\sum \nu'(l_i)\otimes
l^i}\right)/[X]-\frac{1}{8}\sigma(X)=
$$
$$=-\frac{1}{8}\sigma(X)+\frac{1}{8}c^2+\frac{1}{3!}(\frac{c}{2}+ \sum
\nu'(l_i)\otimes l^i)^3/[X]+
\frac{1}{4!}(\frac{c}{2}+\sum \nu'(l_i)\otimes l^i)^4/[X]  \ .
$$
To simplify this expression, put $c_{ij}:=\langle c\cup l^i\cup l^j,[X]\rangle$,
$l_{hijk}:=\langle l^h\cup l^i\cup l^j\cup l^k,[X]\rangle$. The numbers
$c_{ij}$ are
even  since $c$ is characteristic and $(l^i\cup l^j)^2=0$. Substituting into the
formula above we find:
$$\begin{array}{cl}
ch(index(\Dr_\A))=&-\frac{1}{8}\sigma(X)+\frac{1}{8}c^2+\frac{1}{2}
\sum\limits_{i<j}c_{ij}\ [\nu'(l_i)\cup
\nu'(l_j)]+\\ \\
&+\sum\limits_{h<i<j<k} l_{hijk}[\nu'(l_h)\cup\nu'(l_i)\cup \nu'(l_j)\cup
\nu'(l_k)]\ .
\end{array}$$
The cohomology class
$$
u_c:=\frac{1}{2}  \sum_{i<j}c_{ij}\
[\nu'(l_i)\cup \nu'(l_j)]\in H^2({\cal T}'_{\beta}(L),\Z) $$
has the following invariant description: The assignment
$$
(a,b)\longmapsto \ha\langle c\cup a\cup b,[X]\rangle
 $$
defines a $\Z$-valued skew-symmetric  bilinear  form  on
$H^1(X,\Z)$; using the isomorphism $H^1(X,\Z) {\simeq} 2H^1(X,\Z)$, we get a
 cohomology  class in  $H^2({\cal
T}'_{\beta}(L),\Z)$ $=\Lambda^2(2H^1(X,\Z)^{\vee})$, and this coincides with
$u_c$. Clearly $u_c=c_1(index(\Dr_{\A}))$.
\begin{lm} { } [LL] Let $X$ be a 4-manifold with with $b_+=1$. Then
$$c_k(index(\Dr_{\A}))=\frac{1}{k!}u_c^k \ .
$$
\end{lm}
\pf  In the case $b_+=1$, all coefficients $l_{hijk} $ vanish, so
$ch_k(index(\Dr_{\A}))=0$ for $k\geq 2$.
\qed

Regard now $u_c\in H^2({\cal T}'_{\beta}(L),\Z)$ as an element of
$\Lambda^2\left(\qmod{H_1(X,\Z)}{\rm Tors}\right)$.
\begin{pr} Let $X$ be a closed connected oriented 4-manifold with $b_+=1$, and
$\oo_1$ an orientation of $H^1(X,\R)$. Let $c\in H^2(X,\Z)$ be a characteristic
element, $\hat P$ a compatible $Spin^c(4)$-bundle and
$\cg\in\pi_0\left(\qmod{{\cal C}}{\Aut(\hat P)}\right)$. Choose the generator
$l_{\ooo_1}\in\Lambda^{b_1}(H^1(X,\Z))$ which defines the orientation $\oo_1$.
For every $\lambda\in\Lambda^r\left(\qmod{H_1(X,\Z)}{\rm Tors}\right)$ with
$r\equiv b_1$ (mod 2) and $0\leq r\leq\min(b_1,w_c)$, we have:
$$\sigma_{X,\ooo_1}(\cg)(\lambda)= \left\langle\lambda\wedge\exp(-u_c)
 ,l_{\ooo_1}\right\rangle \ .
 $$
 In all other cases
$\sigma_{X,\ooo_1}(\cg)(\lambda)=0$.
\end{pr}

The proof follows from a more general result which we will now explain. \\

Let $T$ be a closed  connected oriented  manifold. Consider Hilbert vector
bundles $E$ and $F$  over $T$, and a   smooth family of Fredholm operators
$q_t:E_t\map  F_t$ of constant  index $\delta$.    Suppose  the   map
$\tilde q: E\setminus\{0\} \map  F$ is a submersion, so that its zero locus
$\tilde
T:=Z(\tilde q)$ is a (finite dimensional)  manifold which fibers over the
possibly
singular  "Brill-Noether locus"
$BN_q:=\{ t\in T|\ \ker q_t\ne \{0\}\}$ of the family $q$.    Put $\hat
T:=\qmod{\tilde
T}{\C^*}$. The projection $\hat p:\hat T\map  T$ induces a projective
fibration over
$BN_q$. Note that
$\hat T$ comes with a canonical cohomology class $u\in H^2(\hat T,\Z)$
induced by the
dual of the  $\C^*$-bundle  $\tilde T\map \hat T$; the restriction of $u$
to  any fiber
$\hat p^{-1}(t)=\P(\ker q_t)$ is the   positive generator of the  second
cohomology group
of this projective space.    We wish to calculate the direct images
$\hat p_*(u^k)$ for all $k\in\N$, $k\geq\delta$ in terms of the Chern
classes of the (virtual)
index bundle $index(q)$ of the family $q$. In the particular case of  Dirac
operators on a
4-manifold with $b_+=1$,   similar computations have been carried out in
[LL] .
\begin{pr}Let $c_i=c_i(index(q))$ be the Chern classes of   $index(q)$, and
define polynomials
 $(p_k)_{k\geq \delta-1}$ by the recursive relations:
$$p_{\delta-1}=1,\ p_k=-\sum_{i=1}^{k-\delta+1} c_i p_{k-i} \ .
$$
For every  non-negative integer $k\geq \delta$ we have
$$\hat p_*(u^k)= p_k(c_1,c_2,\dots) \ ,
$$
hence
$\hat p_*(u^{\delta-1})= 1\in H^0(T,\Z)$ when $\delta>0$.
\end{pr}

\pf One can  find a smooth family of Fredholm operators  $(Q_t)_{t\in T}$,
$Q_t:E_t\oplus
\C^n\map F_t $ with $Q_t|_{E_t\times\{0\}}=q_t$, such that $Q_t$ is
surjective of
positive index $n+\delta$    for every $t\in T$. The associated map
$\tilde Q : (E\oplus\C^n)\setminus\{0\} \map F$ is  a submersion
and $Z(\tilde Q)$  is a locally trivial fiber bundle over $T$ with standard
fiber $\C^{n+\delta}\setminus\{0\}$. Indeed, $Z(\tilde Q)$  is the
complement of the
zero section of the vector bundle
$$V:=\union\limits_{t\in T}\ker Q_t\ .$$
 The space $\tilde T$ can be identified with  the
zero locus $Z(\zeta)$, of the map
$$\zeta:Z(\tilde Q) \map \C^n$$
given by  $\zeta(e,z)=z$. The map $\zeta$ is a submersion in all points of
$\tilde T$, since $\tilde q$ was such.

Let $p_V:\P(V)\map T$ be the obvious projection, and denote by $U\in
H^2(\P(V),\Z)$ the Chern
class of the dual  of the tautological bundle. The map $\hat p:\hat T\map
T$  factors through
the inclusion $j:\hat T\hookrightarrow\P(V)$, and the fundamental class
$j_*[\hat T]$
is   Poincar\' e dual to $U^n$. Therefore  we have
$$\hat p_*(u^k)=[PD_T^{-1}\circ  p_{V*}\circ j_*\circ PD_{\hat
T}](u^k)=[PD_T^{-1}\circ  p_{V*}\circ PD_{\P(V)}](U^{k+n})= p_{V*}(U^{k+n}) .$$
Since $index(q)=[V]-[\C^n]\in K(T)$, we have $c_i(V)=c_i(index(q))=c_i$,
and therefore
$$U^{\delta+n}= -\sum_{i=1}^{\delta+n}
p_{V}^*(c_i)U^{\delta+n-i}  \ .
$$
Multiplying with $U^{k-\delta}$ and using $ p_{V*}(U^{\delta+n-1})=1$, we
get  the recursion relations
$$  p_{V*}(U^{k+n})=-\sum\limits_{i=1}^{\delta+n}p_{V*}(U^{k+n-i}) c_i\ ,$$
hence $\hat p_{V*}(U^{k+n})=p_k$ for  $k\geq \delta-1$ by induction.
\qed

Now we can prove  Proposition 14  by applying  the result above to the map
$\hat p\hat{\cal T}_\beta'^\gamma(L)\map{\cal T}'_\beta(L)$ and the family
$\Dr_{\A}$
of Dirac operators over
${\cal T}'_\beta(L)$.

Since $c_k(index(\Dr_{\A}))=\frac{1}{k!} u_c^k$, we get $ p_*(u^{\delta-1+k})=
p_{\delta-1+k}=\frac{(-1)^k}{k!} u_c^k$, hence
$$\hat p _*(u^{\frac{w_c-r}{2}})=
\frac{ (-1) ^{\left[\frac{b_1-r}{2}\right]} }{\ \
\left[\frac{b_1-r}{2}\right]!\ }\ u_c\
^{\left[\frac{b_1-r}{2}\right]}
$$
for   any  non-negative integer $r$ with $ r\leq
\min(b_1,w_c)$ and $r\equiv b_1$ (mod 2).

Therefore, for every $\lambda\in \Lambda^r\left(\qmod{H_1(X,\Z)}{\rm
Tors}\right)$, we
find
$$\begin{array}{cl}\sigma_{X,\ooo_1}(\cg)(\lambda)&=\langle \hat
p^*\nu'(\lambda)\cup
u^{\frac{w_c-r}{2}},[\hat {\cal T}_\beta'^\gamma(L)]_{\ooo_1}\rangle=
\langle  \nu'(\lambda)\cup  \hat p_*(u^{\frac{w_c-r}{2}}),[ {\cal
T}_\beta'  (L)]_{\ooo_1}\rangle\\ &= \frac{ (-1)^{\left[\frac{b_1-r}{2}\right]} }{\ \
\left[\frac{b_1-r}{2}\right]!\ }\left\langle  \lambda\wedge u_c\
^{\left[\frac{b_1-r}{2}\right]},l_{\ooo_1}\right\rangle \ ,
\end{array}$$
which proves the proposition.
\qed

\subsection{Wall crossing}

The following theorem generalizes results of [W], [KM] and [LL].

\begin{thry} Let $X$ be a closed connected oriented 4-manifold with $b_+=1$, and
$\oo_1$  an  orientation of $H^1(X,\R)$. For every class $\cg$ of
$Spin^c(4)$-structures of Chern class $c$ and  every component ${\bf H}_0$
of ${\bf
H}$, the following holds:
$$ SW_{X,(\ooo_1,{\bf H}_0)}(\cg)(+) - SW_{X,(\ooo_1,{\bf
H}_0)}(\cg)(-) =\sigma_{X,\ooo_1}(\cg)\ .
$$
\end{thry}

We need some preparations before we can prove the theorem.

Fix a  compatible $Spin^c(4)$-bundle $\hat P$ with spinor bundles $\Sigma^{\pm}$
and determinant $L:=\det \Sigma^{\pm}$. Choose a  Clifford map
$\gamma$ with $\cg_\gamma=\cg$, set  $g=g_\gamma$, and let
 $\omega_g$ be  the  generator of $\H_{+,g}^2(X)$ whose class belongs to ${\bf
H}_0$. Put $s_c:=  c\cdot [\omega_g]$, so that the $g$-harmonic
representant of $c-s_c[\omega_g]$ is $g$-anti-selfdual.

Consider first a cohomology class $b_0\in H^2_{\rm DR}(X)$ with
$(c- b_0)\cdot[\omega_g]=0$, and let $\beta_0 \in b_0$. The moduli space
${\cal W}_{\beta_0}:={\cal W}_{X,\beta_0}^\gamma$ contains the closed subset
 of reducible solutions of the form
$(A,0)$, where $A$ solves the   equation
$$F_A^++2\pi i\beta_0^+=0 \ .
$$
This closed subspace can be identified with ${\cal
T}_{ \beta_0^+}'(L)$.

 Consider  the equations
$$\left\{\begin{array}{lll}\Dr_A\Psi&=&0\\
\Gamma(F_A^++2\pi i\beta^+)&=&2(\Psi\bar\Psi)_0
 \end{array}\right. \eqno{( {SW}^\gamma)}
 $$
for a triple $(A,\Psi,\beta)\in {\cal A}(L)\times
A^0(\Sigma^+)\times Z^2_{\rm DR}(X)$,  denote     by ${\cal W}$  the
corresponding moduli space of solutions, and by $p:{\cal W}\map  Z^2_{\rm
DR}(X)$
the natural projection.
${\cal W}$ is  singular in the  points of the form
$[A,0,\beta]$.  If such a triple is a solution, then
$(c-[\beta])\cdot[\omega_g] =0$,
and the singular part of ${\cal W}$ is
$${\cal S}=\union\limits_{(c-[\beta])\cdot[\omega_g]}{\cal T}_{ \beta^+}'(L) \ .
$$
Now perform a "real blow up"  of the singular locus ${\cal S}\subset{\cal
W}$  in
the $\Psi$-direction. This means, consider the equations
$$\left\{\begin{array}{lll}\Dr_A\Phi&=&0\\
\Gamma(F_A^++2\pi i\beta^+)&=&2t(\Phi\bar\Phi)_0
 \end{array}\right. \eqno{(\hat{SW}^\gamma)}
 $$
for a tuple $(A,\Phi,t,\beta)\in {\cal A}(L)\times
S(A^0(\Sigma^+))\times\R\times Z^2_{\rm DR}(X)$, where $S(A^0(\Sigma^+))$ is
the unit sphere in $A^0(\Sigma^+)$ with respect to the $L^2$-norm. Denote  by
$\hat{\cal W}$ the moduli space of solutions of
$(\hat{SW}^\gamma)$ and by
$q_\R$, $q$ the natural projections on $\R$ and   $Z^2_{\rm DR}(X)$
respectively.
Let
$\hat {\cal W}^{\geq 0}:=q_\R^{-1}(\R_{\geq 0})$ be the closed subset defined by
the inequality $t\geq 0$, and, for a  form $\beta\in Z^2_{\rm DR}(X)$, put
$\hat{\cal W}_{\beta}:=q^{-1}(\beta)$ and $\hat{\cal W}_{\beta}^{\geq
0}=\hat{\cal
W}_{\beta}\cap \hat {\cal W}^{\geq 0}$.

There is a natural map  $\hat {\cal W}^{\geq 0}\textmap{\rho}{\cal W}$,
given by
$(A,\Phi,t,\beta)\longmapsto (A,t^{\frac{1}{2}}\Phi,\beta)$, which contracts the
locus $\hat{\cal S}:=\{t=0\}$ to ${\cal S}$ and defines a real analytic
isomorphism
$\hat {\cal W}^{\geq 0}\setminus \hat{\cal S}\map {\cal W} \setminus{\cal S}$.
Note  that
$$\hat{\cal S}=\union\limits_{(c-[\beta])\cdot[\omega_g]}\hat{\cal T}_{\beta
^+}'(L)\ ,
$$
where $\hat{\cal T}_{\beta ^+}'(L)$ is a projective fibration over the
Brill-Noether locus in ${\cal T}_{\beta ^+}'(L)$.
\begin{lm} If   $(g,[\beta])$ is $c$-good, then $q_\R|_{\hat{\cal
W}_{\beta}^{\geq 0}}$ is bounded below by a
positive number. The space $\hat{\cal
W}_{\beta}^{\geq 0}$ is  open and closed in $\hat{\cal W}_{\beta}$, and it is
isomorphic to ${\cal W}_\beta$ via the map $\rho$.
\end{lm}
\pf If $\hat{\cal W}_{\beta}$ would contain a
sequence $[(A_n,\Phi_n,t_n,\beta)]_{n\in\N}$  with
$t_n\searrow 0$ , then, by Proposition 1 $i)$,  we could find a  subsequence
$(m_n)_{n\in\N}\subset\N$ such that $(A_n,t_n^{\frac{1}{2}}\Phi_n)_{n\in\N}$
converges to a  point in
${\cal W}_\beta$. This point must be reducible.

\qed
\begin{lm}   $\hat{\cal W}$    is a smooth manifold .
\end{lm}
\pf Let $(\Psi^-,S)$ be orthogonal to $\im d_\tau(\hat F)$, where $\hat
F=(\hat F_1,\hat F_2)$ is  the map given by the left-hand side of $\hat
{(SW^\gamma)}$, and $\hat\tau=(A,\Phi,t,\beta)$ is a point in ${\cal
A}(L)\times S(A^0(\Sigma^+))\times\R\times Z^2_{\rm DR}(X)$. Using
 variations of
$\beta$, we get  $i\Gamma^{-1}(S)\in  d^*(A^3)\cap A^2_+$, hence $S=0$.
Using now variations of $A$ and the non-triviality of $\Phi$,  we get
$\Psi^-=0$.
\qed
\begin{lm}
The linear map $h:A^1\times\R\map Z^2_{\rm DR}(X)$, given by
$h(\alpha,s)=s\omega_g+d\alpha$,  is transverse to $q$.
\end{lm}
\pf  The image of $dh$ is the subspace $d(A^1)\oplus\H^2_{+,g}(X)$ of $Z^2_{\rm
DR}(X)$ which coincides with the orthogonal complement  of $\H^2_{-,g}(X)$ in
$Z^2_{\rm DR}(X)$. It suffices to show that $\H^2_{-,g}(X)$ is contained in the
image of $d_{\hat\tau}q$, for every $\hat\tau\in\hat{\cal W}$. But if
$\hat\tau=(A,\Phi,t,\beta)$ solves the equations
$(\hat {SW}^\gamma)$, then also
$\hat\tau_t:=(A,\Phi,t,\beta+r\omega_-)$ is   a solution for every
$\omega_-\in\H_{-,g}(X)$ and every $r\in\R$. Therefore $\H_{-,g}(X)$ is
contained
in the image of $d_{\hat\tau} q$.
\qed

Since $h$ is transverse to $q$, the fibre product ${\cal
V}:=(A^1\times\R)\times_{(h,q)}
\hat{\cal W}$ is a smooth manifold, and, putting
$h_\alpha:=h(\alpha,\cdot)$ , we see that for any
$\alpha$ in  a    second category subset of $A^1$, the fibre product  ${\cal
V}_\alpha:=\R\times_{(h_\alpha,p)}\hat{\cal W}$ is   a smooth submanifold of
${\cal V}$.
\begin{lm} The map $\theta:{\cal V}\map \R$ , projecting
$(\alpha,s,A,\Phi,t,s\omega_g+d\alpha)$ to $t$, is a submersion.   For any
$\alpha$ in a   second category subset of $A^1$, the restricted map
$\theta|_{{\cal V}_\alpha}$ is a submersion in all points of $Z(\theta)\cap{\cal
V}_{\alpha}$.
\end{lm}

\pf  ${\cal V}$ can be identified with the moduli space of tuples
$(\alpha,s,A,\Phi,t)\in A^1\times\R\times  {\cal
A}(L)\times S(A^0(\Sigma^+))\times\R$ satisfying the equations
$$\left\{\begin{array}{ccl}  \Dr_A\Phi&=&0\\
\Gamma(F_A^++2\pi i( s\omega_g+d^+\alpha))-2t(\Phi\bar\Phi)_0&=&0 \ .\\
\end{array}\right.
$$
Since the map defined by the left hand side of this system is a submersion
in every
tuple solving the equations,  it suffices to show that the map
$$T:A^1\times\R\times  {\cal
A}(L)\times S(A^0(\Sigma^+))\times\R\map  A^0(\Sigma^-)\times
A^0(su(\Sigma^+))\times\R\ ,$$
 defined by
$$T(\alpha,s,A,\Phi,t)=\left(\begin{array}{c}\Dr_A\Phi\\
\Gamma(F_A^+ + 2\pi i( s\omega_g+d^+\alpha)) -2t(\Phi\bar\Phi)_0 \\ t
\end{array}\right)\ ,
$$
 is a submersion in every point $v=(\alpha,s,A,\Phi,t)$ with
$[(\alpha,s,A,\Phi,t)]\in{\cal V}$. If
$(\Psi^-,S,r)$ is orthogonal to
$\im(d_v T)$, use first variations of $s$ and $\alpha$ to get $S=0$,   then
variations of $A$ to get
$\Psi^-=0$, and then variations of $t$ to get $r=0$.

The second assertion follows by applying Sard's Theorem to the projection
$Z(\theta)\map  A^1$   onto the first factor.
\qed
\begin{lm} For any $\alpha$ in a second category subset of $A^1$, the
moduli space
$\hat {\cal T}_{(s_c\omega+d^+\alpha)}'(L)$ is a smooth manifold.
\end{lm}

\pf If $(\alpha,s,A,\Phi,0)\in Z(\theta)$, then the pair
$(g,[s\omega_g+d\alpha])$
cannot be $c$-good, hence the $s$-component of every point in $Z(\theta)$ must
be $s_c$. Thus  there is a natural identification
$\hat {\cal T}_{(s_c\omega+d^+\alpha)}'(L)=Z(\theta)\cap{\cal
V}_{\alpha}$.
\qed

Since a  pair  $([\omega_g],[h(\alpha,s)])$ belongs to the wall $c^{\bot}$
iff $s=s_c$, every point   $(\alpha,s,A,\Phi,t,s\omega_g+d\alpha)\in {\cal
V}$ with $s\ne s_c$  must have a non-vanishing $t$-component.
\begin{lm} The map $\chi:{\cal V}\map \R$,  projecting
$(\alpha,s,A,\Phi,t,s\omega_g+d\alpha)$ to $ s$, is a submersion in every
point  of  ${\cal V}\setminus Z(\theta)$, in particular in all points
$(\alpha,s,A,\Phi,t,s\omega_g+d\alpha)$ with $s\ne s_c$.
\end{lm}

\pf It suffices to show that the map
$$U:A^1\times\R\times  {\cal
A}(L)\times S(A^0(\Sigma^+))\times\R\map  A^0(\Sigma^-)\times
A^0(su(\Sigma^+))\times\R\ , $$
 defined by
$$U(\alpha,s,A,\Phi,t)=\left(\begin{array}{c} \Dr_A\Phi\\
\Gamma(F_A^++2\pi i( s\omega_g+d^+\alpha)) -2t(\Phi\bar\Phi)_0 \\ s
\end{array}\right)\ ,
$$
 is a submersion in every point $v=(\alpha,s,A,\Phi,t)$ with $t\ne 0$ . If
$(\Psi^-,S,r)$ is orthogonal to
$\im(d_v U)$, we first use first variations of  $\alpha$ to see that $S$ is
orthogonal to $d^+(A^1)$, then variations of $A$ to get
$\Psi^-=0$, and then
variations of $\Phi$ and $t$  ( $t\ne 0$ !) to get $S=0$.    Finally, using
variations of
$s$ one finds $r=0$.
\qed

Applying Sard's theorem again, we have
\begin{lm} Let $S_0$ be a countable  subset of $\R\setminus\{s_c\}$. Then,
for every
$\alpha$ in a   second category subset of
$A^1$, the restricted map
$\chi|_{{\cal V}_\alpha}$ is a
submersion in every point of $\union\limits_{s\in S_0}Z(\chi|_{{\cal
V}_\alpha}-s)$
\end{lm}
\vspace{3mm}

Now we can prove the theorem. \\

\pf Fix $\alpha\in A^1$ such that $\theta|_{{\cal V}_\alpha}$ is a
submersion in every
point of $Z(\theta|_{{\cal V}_\alpha})$ and $\chi|_{{\cal V}_\alpha}$ is a
submersion in
every point of  $Z(\chi|_{{\cal V}_\alpha}-(s_c\pm 1))$. Set
$\beta_0=h_\alpha(s_c)=s_c\omega_g+d\alpha$,
$b_0:=[\beta_0]=[s_c\omega_g]$, and $\beta_\pm:=h(\alpha, s_c\pm 1)$,
$b_\pm:=[\beta_\pm]$. Then
$([\omega_g],b_0)$  belongs to the wall  $c^{\bot}$, and the intersections $\hat
{\cal W}_{\beta_\pm}:=\hat{\cal W}\cap\{\beta=\beta_\pm\}$ are smooth. Note
 that $([\omega_g],b_{\pm})\in C_{{\bf H}_0,\pm}$.

The space
$$\bar{\cal V}^{\geq 0}:={\cal V}_\alpha  \cap\{s_c-1\leq
s\leq s_c+1\}\cap\{t\geq 0 \} = [s_c-1,s_c+1]\times_{(h_\alpha,p)}\hat {\cal
W}^{\geq 0}$$
 is a
smooth manifold with boundary  $\hat {\cal W}^{\geq 0}_{\beta_+}\union
\hat {\cal W}^{\geq 0}_{\beta_-}\union \hat{\cal T}_{\beta_0^+}'(L)$ which
is isomorphic to  ${\cal W}_{\beta_+}\union {\cal W} _{\beta_-}\union \hat{\cal
T}_{\beta_0^+}'(L)$ according to Lemma 17.
\dfigure 181mm by 177mm (bordy scaled 500 offset 0mm:)

Indeed, by the choice of $\alpha$,  ${\cal V}_\alpha$ is smooth,
the projection on the
$t$ component is a submersion in all points of $\hat{\cal
T}_{\beta_0^+}'(L)$ (Lemma 20), and   the projection on the $s$-component
is a submersion in the points of
$\hat {\cal W}^{\geq 0}_{\beta_\pm}$ (Lemma 22) .

Now use the orientation $\oo_1$ of $H^1(X,\R)$ and the component ${\bf
H}_0$ of ${\bf H}$ to endow the smooth moduli spaces
$\hat {\cal W}^{\geq 0}_{\beta_\pm}={\cal W} _{\beta_\pm}$
with the corresponding orientations. The manifold
with  boundary  $\bar{\cal V}^{>0}$ can be oriented by  $\oo_1$,  ${\bf
H}_0$, and by choosing the   natural  orientation of the $s$-coordinate.
Then  the
oriented boundary of  $\bar{\cal V}^{>0}$ is
$$\partial\bar{\cal V}^{>0}=\hat {\cal W}^{\geq 0}_{\beta_+}\union
\left(-\hat {\cal W}^{\geq 0}_{\beta_-}\right)\ .$$
Recall that the moduli space $\hat{\cal T}_{\beta_0^+}'(L)$ could be
 oriented using \underbar{only} the orientation $\oo_1$ of $H^1(X,\R)$.

To determine the sign of the part
$\hat{\cal T}_{\beta_0^+}'(L)$ of the oriented boundary $\partial\bar{\cal
V}^{\geq 0}$ is a technical problem, which can be solved by a careful
examination
of the Kuranishi  model  for   ${\cal V}_\alpha$ in a  point  of $\hat{\cal
T}_{\beta_0^+}'(L)$.  The final result is
$$\partial\bar{\cal V}^{\geq 0}=\hat {\cal W}^{\geq 0}_{\beta_+}\union
\left(-\hat {\cal W}^{\geq 0}_{\beta_-}\right)-\hat{\cal
T}_{\beta_0^+}'(L)\ .$$

Since the cohomology classes $u$, $\nu(l_i)$ extend to $\hat
{\cal W}$ and ${\cal V}$, and   their restrictions to the moduli space
$\hat{\cal
T}_{\beta_0^+}'(L )$ coincide with the corresponding classes defined  in the
section above,  the theorem follows from the relation
$$[\hat {\cal W}^{\geq 0}_{\beta_+}]_{\ooo}-[\hat {\cal W}^{\geq
0}_{\beta_-}]_{\ooo}-[\hat{\cal T}_{\beta_0^+}'(L)]_{\ooo_1}=0
$$
between the fundamental classes of the three moduli spaces. Here $\oo$ is the
orientation determined by $\oo_1$ and ${\bf H}_0$.

\qed\\
{\bf Remark:} Let $(X,g)$ be a manifold of positive scalar curvature,
and let $c\in H^2(X,\Z)$ be a characteristic element such that $(g,0)$ is
$c$-good.
Then
$SW_{X,(\oo_1,{\bf H}_0)}(\cg)(\cdot)=0$ for at least one element in
$\{\pm\}$.
This element is determined by the sign in the inequality $\pm
c\cdot[\omega_g]<0$,
and    the other value of
$SW_{X,(\oo_1,{\bf H}_0)}(\cg)(\cdot)$ is determined by the wall crossing
formula.

\section { Seiberg-Witten invariants of K\"ahler surfaces}
 Let $(X,g)$ be a
K\"ahler surface with K\"ahler form
$\omega_g$, and let $\cg_0$ be the class of the canonical
$Spin^c(4)$-structure of
determinant
$K_X^{\vee}$ on $(X,g)$. The corresponding spinor bundles are
$\Sigma^{+}=\Lambda^{00}\oplus\Lambda^{02}$, $\Sigma^{-}=\Lambda^{01}$  [OT1].
There is a natural bijection between classes of $Spin^c(4)$-structures $\cg$ of
Chern class $c$ and isomorphism classes of line bundles
$M$ with $2c_1(M)-c_1(K_X)=c$. We denote by $\cg_M$ the class defined by a
line bundle  $M$. The spinor bundles of  $\cg_M$  are the tensor products
$\Sigma^{\pm}\otimes M$, and the map $\gamma_M:\Lambda^1_X\map \R
SU(\Sigma^+\otimes M,\Sigma^-\otimes M)$ given by
$\gamma_M(\cdot)=\gamma_0(\cdot)\otimes\id_M$ is a Clifford map representing
$\cg_M$.

Let $C_0$ be the Chern connection in the anti-canonical bundle $K_X^{\vee}$. We
 use the variable substitutions $A:=C_0\otimes B^{\otimes 2}$ with $B\in{\cal
A}(M)$ and $\Psi=:\varphi+\alpha\in A^0(M)\oplus A^{02}(M)$ to rewrite the
Seiberg-Witten equations for $(A,\Psi)$ in terms of $(B,\varphi+\alpha)\in{\cal
A}(M)\times[A^0(M)\oplus A^{02}(M)]$.
\begin{pr} Let $(X,g)$ be a  K\"ahler surface, and $\beta\in A^{1,1}_\R$ a
closed
real $(1,1)$-form in the class $b$. Let  $M$  be a  Hermitian line bundle such
that $(2c_1(M)-c_1(K_X)-b)\cdot[\omega_g]<0$. A pair
$(B,\varphi+\alpha)\in {\cal A}(M)\times\left[A^0(M)\oplus A^{02}(M)\right]$
solves the $\beta$-twisted Seiberg-Witten equations
$(SW^{\gamma_M}_\beta )$ iff:
$$
\left\{\begin{array}{l}F_B^{20}=F_B^{02}=0\  \\
\alpha=0\ ,\ \ \bar\partial_B(\varphi)=0 \ \\
i\Lambda_gF_B+\frac{1}{2}\varphi\bar\varphi+(\frac{s}{2}-\pi\Lambda_g\beta)=0
\ .\end{array}\
\right. $$
\end{pr}
\pf  The pair $(B,\varphi+\alpha)$ solves $(SW^{\gamma_M}_\beta )$ iff
$$\begin{array}{ll}F_A^{20}&=-\varphi\otimes\bar\alpha\\
F_A^{02}&=\ \alpha\otimes\bar\varphi\\
\bar\partial_B(\varphi)&=\ i\Lambda\partial_B(\alpha) \\
i\Lambda_g(F_A+2\pi i\beta)&
=- \left(\varphi\bar\varphi-*(\alpha\wedge\bar\alpha)\right).\end{array}\
$$

Using Witten's transformation $(B,\varphi+\alpha)\longmapsto(B,\varphi-\alpha)$,
we find
$\varphi\otimes\bar\alpha=\alpha\otimes\bar\varphi=0$, hence
$F_A^{20}=F_A^{02}=0$, so that
$\varphi$ or $\alpha$ must vanish. Putting $c:=2c_1(M)-c_1(K_X)$ and integrating
the last equation over $X$ we get
$$\frac{1}{2\pi}\int\limits_X(|\alpha|^2-|\varphi|^2)\frac{\omega_g^2}{2} =
\int\limits_X(\frac{i}{2\pi}F_A-\beta)\wedge\omega_g=(c-b)\cup[\omega_g]<0
\ ,$$
hence  $\alpha=0$.
\qed

 Let  ${\cal D}ou(m)$ be the Douady space of all effective
divisors $D$ on $X$ with $c_1({\cal O}_X(D))=m$.
\begin{thry}
Let $(X,g)$ be a  connected K\"ahler surface,   and let $\cg_M$ be the class of
the $Spin^c(4)$-structure associated to a  Hermitian line bundle
$M$ with $c_1(M)=m$. Let $\beta\in A^{1,1}_{\R}$ be a   closed  form
representing
the class $b$ such that
$ (2m-c_1(K_X)-b)\cup[\omega_g]<0$ ($>0$).
\hfill{\break}
i) \ If $c\ \not\in\  NS(X)$, then ${\cal W}_{X,\beta}^{\gamma_M}
=\emptyset$. If $c\in NS(X)$, then there is a natural real analytic
isomorphism
${\cal W}_{X,\beta}^{\gamma_M}\simeq  {\cal
D}ou(m)$ $({\cal D}ou(c_1(K_X)-m))$.
 \hfill{\break}
ii)  ${\cal W}_{X,\beta}^{\gamma_M} $ is  smooth at a point corresponding to
$D\in{\cal D}ou(m)$ iff
$h^0({\cal O}_D(D))=\dim_D{\cal D}ou(m)\ .$
This condition is always satisfied when $h^1({\cal O}_X)=0$. \hfill{\break}
iii) If  ${\cal W}_{X,\beta}^{\gamma_M}$ is  smooth at a point
corresponding to $D$,
then it has the expected dimension in this point iff $h^1({\cal O}_D(D))=0$.
\end{thry}
\pf Clearly $c\in NS(X)$ is a necessary condition for ${\cal
W}_{X,\beta}^{\gamma_M}\not =\emptyset$. Putting again $c:=2m-c_1(K_X)$, we
may assume that we are in the case $(c-b)\cdot[\omega_g]<0$, since the other one
can be reduced to it by Serre duality. Under this assumption ${\cal
W}^{\gamma_M}_{X,\beta}$ can be identified with the moduli space of holomorphic
pairs $(\bar\partial,\varphi)\in{\cal H}(M)\times A^0(M)$ for which the
generalized
vortex equation
$$i\Lambda_g F_h+\frac{1}{2}\varphi\bar\varphi^h
+(\frac{s}{2}-\pi\Lambda_g\beta)=0
$$
is solvable. The latter space is naturally isomorphic with the Douady space
${\cal
D}ou(m)$ [OT1]. The remaining assertions follow from the long exact cohomology
sequence of the structure sequence
$$0\map {\cal O}_X\textmap{\cdot D}{\cal O}_X(D)\map  {\cal O}_D(D)\map  0\ .
$$
\qed

{\bf Remark:}  We use the complex structure of the surface to orient
$H^1(X,\R)$ and
$\H^2_{+,g}(X)$. With this convention, the following holds:

The natural isomorphism
${\cal W}^{\gamma_M}_{X,\beta}\simeq  {\cal D}ou(m)$   respects the orientation
when
$(2m-c_1(K_X)-b)\cup[\omega_g]<0$.
If $(2m-c_1(K_X)-b)\cup[\omega_g]>0$, then the isomorphism ${\cal
W}^{\gamma_M}_{X,\beta}\simeq  {\cal D}ou(c_1(K_X)-m)$  multiplies the
orientation by the factor $(-1)^{\chi(M)}$.

The pull-back of the hyperplane class  of  ${\cal D}ou(m)$ is precisely $u$
when
$(2m-c_1(K_X)-b)\cup[\omega_g]<0$.
If
$(2m-c_1(K_X)-b)\cup[\omega_g]>0$, then the pull-back of the hyperplane class
of  ${\cal D}ou(c_1(K_X)-m)$ is $-u$.
\\

Recall that an effective divisor $D$ on a connected complex surface $X$ is
$k$-connected iff $D_1\cdot D_2\geq k$ for every effective decomposition
$D=D_1+D_2$ [BPV]. Canonical divisors of minimal surfaces of Kodaira dimension
$\kappa=1,   2$ are $(\kappa-1)$-connected.
\begin{lm} Every  connected complex surface with $p_g>0$ is oriented
diffeomorphic to a surface which possesses a 0-connected canonical divisor.
\end{lm}
\pf Let $X$ be a connected complex surface with $p_g>0$ and minimal model
$X_{\rm min}$. Let $K_{\rm min}$ be a 0-connected canonical divisor of $X_{\rm
min}$. Choose $b_2(X)-b_2(X_{\rm min})$ distinct  points $x_i\in X_{\rm
min}\setminus{\rm supp}(K_{\rm min})$, and note that $X$ is diffeomorphic to the
blow  up $\hat X_{\rm min}$ of  $X_{\rm min}$ in these points. Denote by
$\sigma$
the projection $\sigma:\hat X_{\rm min}\map X_{\rm min}$  and by $E$ the
exceptional divisor. Then $\hat K:=\sigma^*(K_{\rm min})+E$ is a canonical
divisor
on $\hat X_{\rm min}$. If $\hat K$ decomposes as $\hat K=D_1+ D_2$, then every
component of $E$ is contained in precisely one of the summands, and $K_{\rm
min}=\sigma_*(\hat K)=\sigma_* D_1+\sigma_* D_2$ is a decomposition  of
$K_{\rm min}$. This implies $D_1\cdot D_2=\sigma_* D_1\cdot \sigma_* D_2\geq
0$.
\qed
\begin{co} ([W]) All non-trivial Seiberg-Witten invariants of K\"ahler surfaces
with $p_g>0$ have index 0.
\end{co}
\pf Let $X$ be a K\"ahler surface with $p_g>0$. We may suppose that $X$
possesses a 0-connected canonical divisor $K$, defined by a holomorphic 2-form
$\eta$. Using a moduli space ${\cal W}^\gamma_{X,\eta}$ to calculate the
invariant as in [W], we find an effective decomposition $K=D_1+D_2$. This
implies
$w_c=-D_1\cdot D_2\leq 0$.
\qed
\begin{co} Let $X$ be a K\"ahlerian surface with $p_g=0$ and    $q=0$.  Endow
  $H^1(X,\R)=0$ with the standard orientation  and let ${\bf H}_0$ be the
component of
${\bf H}$ containing  K\"ahler forms . If $m(m-c_1(K_X))\geq 0$, then  we have
$$SW_{X,{\bf H}_0}(\cg_M)(+)=\left\{\begin{array}{ccc} 1&{\rm if} &{\cal
D}ou(m)\ne\emptyset\\
0 &{\rm if} &{\cal D}ou(m)=\emptyset \ ,
\end{array}\right.
$$
and
$$SW_{X,{\bf H}_0}(\cg_M)(-)= \left\{\begin{array}{ccc} 0&{\rm if} &{\cal
D}ou(m)\ne\emptyset\\
-1 &{\rm if} &{\cal D}ou(m)=\emptyset \ .
\end{array}\right.
$$
\end{co}

\pf   Suppose ${\cal D}ou(m)\ne\emptyset$.  Since $p_g=0$ and $q=0$, we
must have
${\cal D}ou(c_1(K_X)-m)=\emptyset$, hence  $SW_{X,{\bf H}_0}(\cg_M)(-)=0$
by Theorem 10.
Now the wall crossing formula implies $SW_{X,{\bf H}_0}(\cg_M)(+)=1$.  If
${\cal D}ou(m)=\emptyset$, then  $SW_{X,{\bf H}_0}(\cg_M)(+)=0$ by Theorem
10, and
$SW_{X,{\bf H}_0}(\cg_M)(-)=-1$ again by the wall crossing formula.
\qed

An interesting formulation is obtained under the additional assumption
${\rm Tors}_2
H^2(X,\Z)=0$. Then a class
$\cg$ of  $Spin^c(4)$-structures is determined by its Chern class $c$ and
$\frac{c_1(K_X)\pm c}{2}$ makes sense. If $c^2\geq c_1(K_X)^2$, then

$$SW_{X,{\bf H}_0}(c)(+)=\left\{\begin{array}{ccc} 1&{\rm if} &{\cal
D}ou(\frac{c_1(K_X)+c}{2})\ne\emptyset\\ \\
0 &{\rm if} &{\cal D}ou(\frac{c_1(K_X)+c}{2})=\emptyset \ ,
\end{array}\right.
$$
and
$$SW_{X,{\bf H}_0}(c)(-)=\left\{\begin{array}{ccc}  -1 &{\rm if} &{\cal D}ou
(\frac{c_1(K_X)-c}{2})\ne\emptyset
\\ \\ 0&{\rm if} &{\cal
D}ou(\frac{c_1(K_X)-c}{2})=\emptyset\ .  \end{array} \right.
$$
\\

{\bf Example:} Let $X=\P^2$,   let $h\in H^2(\P^2,\Z)$ be the class of the ample
generator, and let  ${\bf H}_0$ be the component of ${\bf H}=\{\pm h\}$
which contains $h$.  The classes of $Spin^c(4)$-structures are labelled by odd
integers $c$, and the index corresponding to $c$ is
$w_c=\frac{1}{4}(c^2-9)$. The
chambers of type $c$ contained in ${\bf H}_0\times H^2_{\rm DR}(X)$ are the
half-lines $C_{{\bf H}_0,\pm}=\{ b\in H^2_{\rm DR}(\P^2)|\ \pm(c-b)\cdot h<0\}$.
Using Theorem 10 to
calculate the moduli spaces ${\cal W}_{\P^2,\beta}^{\gamma_{\frac{c-3}{2}}}$ we
find
$${\cal W}_{\P^2,\beta}^{\gamma_{\frac{c-3}{2}}}\simeq\left\{\begin{array}{lll}
|{\cal O}_{\P^2}(\frac{c-3}{2})|&{\rm if}&[\beta]>c\\
|{\cal O}_{\P^2}(\frac{-c-3}{2})|&{\rm if}&[\beta]<c \ .
\end{array}\right.
$$
Taking into account the orientation-conventions, we get by direct verification:
$$SW_{\P^2,{\bf H}_0}(c)(+)=\left\{\begin{array}{lll}
1&{\rm if}& c\geq 3\\
0&{\rm if}&c<3\ ,
\end{array}\right.\  \
SW_{\P^2,{\bf H}_0}(c)(-)=\left\{\begin{array}{lll}
-1&{\rm if}& c\leq-3\\
0&{\rm if}&c>-3\ .
\end{array}\right.$$
\dfigure 187mm by 186mm (plane scaled 350 offset 0mm:)
For every  $c$, the subspace ${\bf H}_0\times\{0\}$ is    contained in
the chamber on which $SW_{\P^2,{\bf H}_0}(c)$ vanishes.\\ \\
{\bf Remark:} Let $X=\hat\P^2$   be the blow-up of $\P^2$ in $r\geq
3$ points,  and fix a non-negative even integer $w$. There exist infinitely many
solutions  $(d;m_1,\dots,m_r)\in\N^{\oplus(r+1)}$ of the equation
$\frac{1}{2}w=\frac{1}{2}d(d+3)-\sum\limits_{i=1}^r\frac{m_i(m_i+1)}{2}\
.$
For every solution $(d;m_1,\dots,m_r)$  let $M$ be the
underlying line bundle of the linear system $|dL-\sum\limits_{i=1}^r  m_i
E_i|$, and set $c:=2 c_1(M)-c_1(K_X)$. Then $w_c=w$  and we have $SW_{X,{\bf
H}_0}(c)(+)=1$,   $SW_{X,{\bf H}_0}(c)(-)=0$  for the component ${\bf H}_0$
containing K\"ahler classes. Hence there exist infinitely many
characteristic classes $c$ with non-trivial Seiberg-Witten invariants  and
prescribed index $w_c=w$.

\parindent0cm
\vspace{0.7cm}
{{\bf References}}\vskip 10pt
{\small

[BPV] Barth, W., Peters, C., Van de Ven, A.: {\it Compact complex surfaces},
Springer Verlag, 1984.

[DK] Donaldson, S.; Kronheimer, P..: {\it The Geometry of
four-manifolds}, Oxford Science Publications  1990.

[KM] Kronheimer, P.; Mrowka, T.: {\it The genus of embedded surfaces in
the projective plane}, Math. Res. Lett. 1,  1994,  pp. 797-808.

[LL] Li, T.; Liu, A.: {\it General wall crossing formula}, Math. Res. Lett.
2,  1995,
pp. 797-810.

[OT1]   Okonek, Ch.; Teleman A.: {\it The Coupled Seiberg-Witten Equations,
Vortices, and Moduli Spaces of Stable Pairs},  Int. J. Math. Vol. 6, No. 6,
1995,
pp. 893-910.

[OT2] Okonek, Ch.; Teleman A.: {\it Quaternionic monopoles}, Comm. Math. Phys. (to
appear).

[OT3] Okonek, Ch.; Teleman A.: {\it Seiberg-Witten invariants for manifolds with
$b_+=1$}, Comptes Rendus   Acad. Sci. Paris (to appear).

[Ta] Taubes,   C.: {\it The Seiberg-Witten and Gromov invariants} Math.
Res. Lett. 2,  1995, pp. 221-238.

 [T] Teleman, A.: {\it  Moduli spaces of monopoles}, Habilitationsschrift,
 Z\"urich 1996  (in preparation).

[W] Witten, E.: {\it Monopoles and four-manifolds}, Math.
Res. Lett. 1,  1994, pp. 769-796.
\vspace{0.2cm}\\
 Ch. Okonek: Mathematisches Institut, Universit\"at Z\"urich,
Winterthurerstr. 190,
\hspace*{2.2cm} CH-8057 Z\"urich \\
A. Teleman: \ Mathematisches Institut, Universit\"at Z\"urich,
Winterthurerstr. 190,
\hspace*{2.4cm}CH-8057 Z\"urich and\\
\hspace*{2.4cm}Faculty of Mathematics, University of Bucharest\\
\hspace*{2.4cm}e-mail: okonek@math.unizh.ch ; teleman@math.unizh.ch}

\end{document}